\documentclass{amsart}
\usepackage{url}
\usepackage{hyperref}
\usepackage[procnames]{listings}
\usepackage{cite}
\usepackage{color}
\usepackage{ amssymb }
\usepackage{amsmath}
\usepackage{mathtools}
\usepackage{amssymb}
\usepackage[toc,page]{appendix}
\numberwithin{equation}{section}
\title{Complexity of  strong approximation on the sphere }
\author{Naser T. Sardari}
\date{\today}
 \usepackage {enumerate}
\usepackage{tikz-cd}

\usepackage{tikz}

\address{Department of Mathematics, UW-Madison, Madison, WI 53706}

\email{ntalebiz@math.wisc.edu}

	\newtheorem{thm}{Theorem}[section]
	
	\newtheorem{prop}[thm]{Proposition}
	
	\newtheorem{rem}[thm]{Remark}	
	\newtheorem{lem}[thm]{Lemma}
	\newtheorem{conj}[thm]{Conjecture}

	\newtheorem{cor}[thm]{Corollary}
	
	\theoremstyle{defi}

	\newtheorem{q}[thm]{sub-sum problem}
	
	\theoremstyle{pf}

\begin{document}
\maketitle

\begin{abstract}
  By assuming some widely-believed arithmetic conjectures, we show that the task of accepting a number that is representable as a sum of $d\geq2$ squares subjected to given congruence conditions is NP-complete.  On the other hand, we develop and implement a deterministic  polynomial-time  algorithm that  represents a number as a sum of 4 squares with some restricted congruence conditions,  by assuming  a  polynomial-time algorithm for factoring integers and Conjecture~\ref{cc}. As an application,    we develop and implement a deterministic  polynomial-time algorithm for navigating LPS Ramanujan graphs,  under the same assumptions.  

\end{abstract}
\tableofcontents
\section{Introduction}

\subsection{Motivtion}
We begin by defining Ramanujan graphs. Fix $k\geq 3$, and let $G$ be a $k$-regular  connected graph with the adjacency matrix  $A_G.$ It follows that $k$ is an eigenvalue of $A_G$. Let $\lambda_G$ be the maximum of the absolute value of all the other eigenvalues of $A_G$. By the Alon-Boppana Theorem~\cite{Lubotzky1988}, $\lambda_G\geq 2\sqrt{k-1}+o(1),$ where $o(1)$ goes to zero as $|G|\to \infty$.  We say that $G$ is a Ramanujan graph, if  $\lambda_G \leq 2\sqrt{k-1}.$ 

 The first construction of Ramanujan Graphs is due to Lubotzky, Phillips and Sarnak~\cite{Lubotzky1988} and independently by Margulis \cite{Margulis}. We refer the reader to \cite[Chapter 3]{Peter}, where a complete history of the construction of Ramanujan graphs and other extremal properties of them are recorded. The LPS construction has the additional property of being strongly explicit. We say that the $k$-regular graph $G$
is strongly explicit, if there is a polynomial-time algorithm that on inputs $\langle v, i\rangle$ where $v\in G,  1 \leq i \leq k$ outputs the (index of the) $i^{th}$ neighbor of $v$. Note that  the lengths of the algorithm's inputs and outputs are $O(\log |G|)$, and so it runs in time $poly\log(|G|).$ This feature of the LPS Ramanujan graphs is very important in their application  to the deterministic error reduction algorithm \cite{Ajtai}; see also \cite{Hoory} for other applications of Ramanujan graphs in Computer Science.  

The main product of this work is a deterministic polynomial-time algorithm for navigating LPS Ramanujan graphs,  by assuming a polynomial-time algorithm for factoring integers and  an arithmetic conjecture, which we formulate next.

 Let $Q(t_0,t_1):=\frac{N}{4q^2}-(t_0+\frac{a_0}{2q})^2-(t_1+\frac{a_1}{2q})^2$, where $q$ is a prime, $N$, $a_0$, and $a_1$ are integers, where  $N\equiv a_0^2+a_1^2 \mod 4q$ and $\gcd(N,4q)=1.$ Define 
\begin{equation}\label{Afr}
A_{Q,r}:=\{(t_0,t_1)\in\mathbb{Z}^2: Q(t_0,t_1)\in \mathbb{Z}, |(t_0,t_1)|<r, \text{ and } Q(t_0,t_1)\geq 0  \},
\end{equation}
where $r>0$ is some positive real number. 
 
%

 
\begin{conj}\label{cc} 
Let $Q$ and $A_{Q,r}$ be as above. There exists constants $\gamma>0$ and $C_{\gamma}>0$, independent of $Q$ and $r$, such that if $ |A_{Q,r}|> C_{\gamma}(\log N)^{\gamma}$ for some $r>0$, then $Q$ expresses a sum of two squares inside $ A_{Q,r}$. \end{conj}

 We denote the following  assumptions by $(*)$:
  \begin{enumerate}\label{assumptions}
 \item  There exists a  polynomial-time algorithm for factoring integers,
 \item  Conjecture~\ref{cc} holds.
 \end{enumerate}


The LPS construction is the  Cayley graphs of $PGL_2(\mathbb{Z}/q\mathbb{Z})$ or $PSL_2(\mathbb{Z}/q\mathbb{Z})$ with $p+1$ explicit generators for every prime $p$ and integer $q$. We denote them by the LPS Ramanujan graph $X^{p,q}$, and the $p+1$ generators by the LPS generators in this paper.  For simplicity for the rest of this paper as in\cite{Lubotzky1988}, we assume that $q\equiv  1$ mod 4 is also a prime, and  is a quadratic residue mod $p$, where $p\equiv  1$ mod 4 is fixed.  By these assumptions,  $X^{p,q}$ is a Cayley graph  over ${\rm PSL}_2(\mathbb{Z}/q\mathbb{Z})$; see Section~\ref{reduction} for the explicit construction of $X^{p,q}$.  We say $v\in X^{p,q}$ is a diagonal vertex,  if it corresponds to a diagonal matrix in ${\rm PSL}_2(\mathbb{Z}/q\mathbb{Z}).$ By a path from $u_1$ to $u_2$, we mean a sequence of vertices $\langle v_0, \dots,v_h\rangle$, where $v_0=u_1$,  $v_h=u_2$, and $v_i$ is connected to $v_{i+1}$ for every $0\leq i \leq h-1.$
\begin{thm}\label{mainthm} Assume~$(*)$. We develop and implement  a  deterministic  polynomial-time algorithm  in $\log(q)$,   that on inputs $\langle u_1, u_2 \rangle$, where $u_1, u_2\in X^{p,q}$ are diagonal vertices,  outputs  a shortest path $\langle v_0, \dots,v_h\rangle$ from $u_1$ to $u_2$.  Moreover, for every $\alpha \geq0$ we have 
\begin{equation}\label{almostdiam}
h\leq \max(\alpha,  3\log_p(q)+\gamma\log_{p}\log(q)+\log_{p}(C_{\gamma}) +\log_p(89)),
\end{equation}
for all, but  at most $89q^4/p^{(\alpha-1)}$ vertices. In particular, for large enough $q$ the distance of any diagonal vertice from the identity is bounded by 
\begin{equation}\label{diamb}(4/3)\log_{p}|X^{p,q}|+\log_p(89).\end{equation}

\end{thm}
 \begin{rem}
  Our algorithm is the $q$-adic analogue of the Ross and Selinger algorithm~\cite{Selinger}, which navigates $PSU(2)$ with a variant of the LPS generators. In their work, the algorithm terminates in polynomial-time under the first assumption in $(*)$, and some heuristic arithmetic assumptions which are implicit in their work. We formulated Conjecture~\ref{cc}, and proved the algorithm terminates in polynomial-time under $(*)$. Moreover, we give quantitative bounds on the size of the output under $(*)$. In particular, \eqref{almostdiam} implies that the distance between all but a tiny fraction of pairs of diagonal vertices is less than $\log_{p}(|X^{p,q}|)+O(\log\log|X^{p,q}|).$ In order to prove our bounds, we introduce a correspondence between the diagonal vertices of $X^{p,q}$ and the index $q$ sublattices of $\mathbb{Z}^2$. This is novel in our work; see Section~\ref{quant}. 
   \end{rem}
It is known that every pair of vertices of a Ramanujan graph  (not necessarily an LPS Ramanujan graph) are connected by a logarithmic number of edges. More precisely,   for any $x,y \in G$, let $d(x,y)$ be the length of the shortest path between $x$ and $y.$ Define the diameter of $G$ by
$
\text{diam}(G) := \sup_{x,y \in G} d(x,y).
$
It is easy to check that $\text{diam}(G) \geq \log_{k-1}|G|$.  If $G$ is a Ramanujan graph then  
$\text{diam}(G) \leq 2\log_{k-1}|G| +O(1)$; see \cite{Lubotzky1988}. Moreover,  we \cite[Theorem 1.5]{Naser} showed quantitatively that all but a tiny fraction of the  pairs of vertices in  $G$ have a distance less than $\log_{k-1}(|G|)+O(\log \log |G|)$.
 
 Bounding the  diameter of the LPS Ramanujan graph $X^{p,q}$ is closely related to the diophantine properties of quadratic forms in four variables \cite{Naser2}. In particular,  we showed that for every prime $p$ there exists an infinite sequence of integers  $\{q_n\}$, such that 
$\text{diam}(X^{p,q_n}) \geq (4/3) \log_{k-1}|X^{p,q_n}|$; see \cite[Theorem 1.2]{Naser}. This shows that our upper bound in \eqref{diamb} is  optimal. 
In fact, by assuming our conjecture  on the optimal strong approximation for quadratic forms in 4 variables~\cite[Conjecture 1.3]{Naser2}, the diameter of $X^{p,q}$ is asymptotically 
 $(4/3) \log_{k-1}|X^{p,q}|$ as $q\to \infty.$  In our joint work with Rivin~\cite{Rivin}, we gave numerical evidences for this asymptotic.  Our navigation algorithm substantially improves the range of our previous numerical results, and gives stronger evidences  for  \cite[Conjecture 1.3]{Naser2}. 
 
  \begin{rem}
  Sarnak in his letter to Scott Aaronson and Andy Pollington~\cite{Sarnak31}  defined the covering exponent of the LPS generators for navigating $PSU(2)$. He conjectured that the covering exponent is $4/3$; see \cite{Naser2} and \cite{Browning}. In particular, this exponent gives the optimal bound on the size of the output of the Ross and Selinger  algorithm. $\lim_{q\to \infty}\frac{\text{diam}(X^{p,q})}{\log_{p}|X^{p,q}|}$ is the $q$-adic analogue of the covering exponent. In fact, \cite[Conjecture 1.3]{Naser2} generalizes Sarnak's conjecture, and it also implies  $$ \lim_{q\to \infty} \frac{\text{diam}(X^{p,q})}{\log_{p}|X^{p,q}|}= 4/3.$$       \end{rem}
%
%

By assuming $(*)$,  we develop a deterministic  polynomial-time algorithm  that returns a short path between every  pair of vertices of  $X^{p,q}$. This version of the algorithm is not restricted to the diagonal vertices, but it does not necessarily return the shortest possible path; see Remark~\ref{remNp}. 

\begin{thm}\label{decomdiag}\normalfont Assume $(*)$. We develop a  deterministic  polynomial-time algorithm  in $\log(q)$,   that on inputs $\langle u_1, u_2 \rangle$, where $u_1, u_2\in X^{p,q}$,  returns a short path $\langle v_0, \dots,v_h\rangle$ from $u_1$ to $u_2$.  Moreover,  we have 
\begin{equation}\label{deter}
h\leq \frac{16}{3}\log_{k-1}|X^{p,q}|+O(1).
\end{equation}
Furthermore, 
\begin{equation}\label{typical}
h\leq 3\log_{k-1}|X^{p,q}|+ O(\log \log (|X^{p,q}|))
\end{equation}
for all but $O(\log(q)^{-c_1})$ fraction of pairs of vertices, where $c_1>0$, and the implicit constant in the $O$ notations and $c_1$ are independent of $q$.  
%
\end{thm}
 We briefly describe our proof in what follows. By \cite[Lemma 1]{Petit2008}, we express any element of ${\rm PSL}_2(\mathbb{Z}/q\mathbb{Z})$ as a product of a bounded number of LPS generators and four diagonal matrices. This reduces the navigation task to the diagonal case, and so Theorem~\ref{mainthm} implies \eqref{deter}.
 
  For proving~\eqref{typical},  we improve on  Lauter, Petit and Quisquater's  diagonal decomposition algorithm. 
 By~\eqref{almostdiam}, the distance of a typical diagonal element from the identity is less than $\log_p|X^{p,q}|+O(\log_p\log(|X^{p,q}|))$. So,  it suffices to show that all but a tiny fraction of vertices are the product of $O(\log_p\log(|X^{p,q}|))$ number of LPS generators and three typical diagonal matrices. It is elementary to see that at least $10\%$ of the vertices of $X_{p,q}$  are the product of a bounded number of LPS generators and three typical diagonal matrices. By the expansion property of the Ramanujan graphs, the distance of all but a tiny fraction of the vertices is less that  $O(\log_p\log(|X^{p,q}|))$ from any subset containing  more than $10\%$ of vertices. This implies  \eqref{typical}. We give the dull details of our argument in Section~\ref{diagsection}.
%

\begin{rem}\label{remNp} By  Theorem~\ref{reductionthm} and Corollary~\ref{cooor}, it  follows that finding the shortest path between a generic pair of vertices is essentially NP-complete; see Remark~\ref{Npprem} for further discussion. The idea of reducing the navigation task  to the diagonal case is due Petit,  Lauter, and Quisquater  \cite{Petit2008}, which is crucial in both  Ross and Selinger~\cite{Selinger}   and this work.  As a result of this diagonal decomposition, the size of the output path is 3 times the shortest possible path for a typical pair of vertices.  Improving the constant  $3$ to $3-\epsilon$ needs  new ideas, and this would have applications in quantum computing. \end{rem}


\subsection{Reduction to  strong approximation on the sphere}\label{reduction}
 \noindent In \cite[Section 3]{Lubotzky1988}, the authors implicitly reduced the task of finding the shortest possible path between a pair of vertices in $X^{p,q}$ to the task of representing a number as a sum of 4 squares subjected to  given congruence conditions, which is the strong approximation on the 3-sphere.  We explain this reduction in this section.

 We begin by  explicitly describing  $X^{p,q}$.  Let $\mathbb{H}(\mathbb{Z})$ denote the integral Hamiltonian quaternions
$$\mathbb{H}(\mathbb{Z}):= \big\{x_0 +x_1i+x_2 j+x_3 k| x_t\in \mathbb{Z } , 0\leq t \leq 3,   i^2=j^2=k^2=-1\big\},$$
where $ij=-ji=k$, etc. Let $\alpha:= x_0 +x_1i+x_2 j+x_3 k\in \mathbb{H}(\mathbb{Z})$. Denote $\bar{\alpha}:=x_0 -x_1i-x_2j-x_3 k$ and $\text{Norm}(\alpha):= \alpha \bar{\alpha}=x_0^2+x_1^2+x_2^2+x_3^2$.
Let 
\begin{equation}
\label{LPSgen}
 S_p:=\{\alpha \in \mathbb{H}(\mathbb{Z}) : \text{Norm}(\alpha)=p, \text{$x_0 > 0$  is odd and $x_1,x_2,x_3$ are even numbers}   \}.
 \end{equation}
 It follows that $S_p=\{\alpha_1, \bar{\alpha_1}, \dots, \alpha_{(p+1)/2}, \bar{\alpha}_{(p+1)/2}   \}.$
 Let $$\Lambda^{\prime}_p:=\{  \beta \in \mathbb{H}: \text{ Norm}(\beta)=p^{h^{\prime}} \text { and } \beta\equiv 1 \text { mod } 2  \}.$$
 $\Lambda^{\prime}_p$ is closed under multiplication. Let $\Lambda_p$ be the set of classes of  $\Lambda^{\prime}_p$ with the relation   $\beta_1\sim\beta_2$  whenever $\pm p^{t_1}\beta_1=p^{t_2}\beta_2$,  where $t_1, t_2 \in \mathbb{Z}$. Then  $\Lambda^{\prime}_p$ form a group with 
$$[\beta_1][\beta_2]=[\beta_1 \beta_2]  \text{ and  }  [\beta][\bar{\beta}]=[1].$$
By \cite[Corollary 3.2]{Lubotzky1988},  $\Lambda_p$ is free on $[\alpha_1],  \dots, [\alpha_{(p+1)/2}]$. Hence,  the Cayley graph of  $\Lambda_p$ with respect to  LPS generator set $S_p$ is  an infinite   $p+1$-regular tree.  LPS Ramanujan graphs are associated to the quotient of this infinite $p+1$-regular tree by appropriate arithmetic subgroups that we describe in what follows. Let $$\Lambda_p(q):=\{  [\beta] \in\Lambda_p:  \beta=x_0 +x_1i+x_2 j+x_3 k \equiv x_0 \text { mod } 2q  \}.$$
$\Lambda_p(q)$ is a normal subgroup of $\Lambda_p$. By \cite[Proposition  3.3]{Lubotzky1988}, since $q\equiv 1$   mod 4 is a prime number and $q$ is a quadratic residue mod $p$,
$$\Lambda_p/\Lambda_p(q)= {\rm PSL}_2(\mathbb{Z}/q\mathbb{Z}) .$$ 
The above isomorphism is defined by sending  $[\alpha]\in \Lambda_p$, to the following matrix $\tilde{\alpha}$ in ${\rm PSL}_2(\mathbb{Z}/q\mathbb{Z})$:
\begin{equation}\label{cores}
\tilde{\alpha}:=\frac{1}{\sqrt{\text{Norm}(\alpha)}} \begin{bmatrix}
x_0+i x_1 & y+i x_3
\\
-y+ix_3 & x_0 -ix_1
\end{bmatrix},
\end{equation}
where $i$ and $\sqrt{p}$ are  representatives of  square roots of $-1$ and $p$  mod $q$.
This identifies the finite $p+1$-regular graph $\Lambda_p/\Lambda_p(q)$ by the Cayley graph of ${\rm PSL}_2(\mathbb{Z}/q\mathbb{Z})$ with respect to  $\tilde{S}_p$ (the image of ${S}_p$ under the above map) that is the  LPS Ramanujan graph $X^{p,q}$.  For $v\in X^{p,q}$, we denote its associated class in  $\Lambda_p/\Lambda_p(q)$ by $[v].$

Finally, we give a theorem which reduces the navigation task on LPS Ramanujan graphs to an strong approximation problem for the 3-sphere. Since $X^{p,q}$ is a Cayley graph, it suffices to navigate from the identity vertex to any other vertex of $X^{p,q}$. 


%
\begin{thm}[Due to Lubotzky, Phillips and Sarnak]\label{reductionthm} 
 Let $v\in X^{p,q}$, and $a_0+a_1 i+a_2j+a_3k\in [v]$ such that $\gcd(a_0,\dots,a_3,p)=1$. There is a bijection between  non-backtracking  paths $(v_0,\dots,v_h)$  of length $h$ from $v_0=id$ to $v_h=v$  in $X^{p,q}$, and the set of integral solutions to the following diophantine equation 
\begin{equation}\label{redd}
\begin{split}
&x_1^2+x_2^2+x_3^2+x_4^2=N,
\\
&x_l\equiv \lambda a_l \text{ mod } 2q \text{ for  $0\leq l \leq 3$ and some $\lambda \in \mathbb{Z}/2q\mathbb{Z}$,    }
\end{split}
\end{equation} 
   where $N=p^h$.
In particular, the distance between $id$ and $v$ in $X^{p,q}$ is the smallest exponent $h$ such that~\eqref{redd} has an integral solution. 
\end{thm}

By~\cite[Conjecture 1.3]{Naser2}, there exists an integral lift if $p^h\gg_{\epsilon} q^{4+\epsilon}$ and $4$ is the optimal exponent. This conjecture implies that $\text{diam}(X^{p,q})$  is asymptotically, 
 $$4/3 \log_{k-1}|X^{p,q}|.$$

\subsection{Complexity of  strong approximation on the sphere} 
 In this section, we give our main results regarding  the complexity of representing a number  as a sum of $d$ squares subjected to given congruence conditions. First, we give our result for $d=2.$

%

\begin{thm}\label{nptheorem}
The problem of accepting $(N,q,a_0,a_1)$ such that the diophantine equation 
\begin{equation*}
\begin{split}
&x_0^2+x_1^2=N,
\\
&x_0\equiv a_0 \text{ and } x_1\equiv a_1 \text{ mod } q,
\end{split}
\end{equation*}
has integral solution $(x_0,x_1)\in \mathbb{Z}^2$ is NP-complete,  by assuming GRH and Cramer's conjecture, or unconditionally by a randomized reduction algorithm.
\end{thm}
The above theorem is inspired by a private communication with Sarnak. He showed us  that the problem of representing a number as a sum of two squares subjected to inequalities on the coordinates is NP-complete, under a randomized reduction algorithm. The details of this theorem appeared in his joint work with Parzanchevski \cite[Theorem 2.2]{Sarnak3}. 

By induction on $d$, we generalize our theorem for every $d\geq 2$.
\begin{cor}\label{cooor}
Let $d\geq 2$. The problem of accepting $(N,q,a_0,\dots,a_{d-1})$ such that the diophantine equation 
\begin{equation}\label{deq}
\begin{split}
&x_0^2+\dots +x_{d-1}^2=N,
\\
&X_0\equiv a_0 \dots x_{d-1}\equiv a_{d-1} \text{ mod } q,
\end{split}
\end{equation}
has integral solution $(x_0,\dots,x_{d-1})\in \mathbb{Z}^d$ is NP-complete,  by assuming GRH and Cramer's conjecture, or unconditionally by a randomized reduction algorithm.
\end{cor}
 On the other hand, by assuming $(*)$ and two coordinates of the congruence conditions  in \eqref{deq} are  zero, we develop and implement a polynomial-time algorithm  for this task for $d=4.$
%

\begin{thm} \label{diagliftt}
Let $q$ be a prime, and  $(a_0,a_1)\in \big(\mathbb{Z}/2q\mathbb{Z}\big)^2$, where $a_0$ is odd and $a_1$ is even.   Suppose that $N=O(q^A),~  \gcd(N,4q)=1$, and  
$a_0^2 +a_1^2 \equiv N  \text{  mod } 4q.$ By assuming $(*)$, we develop and implement a deterministic   polynomial-time algorithm in $\log(q)$ that  finds an integral solution $(x_0,\dots,x_3)\in \mathbb{Z}^4$ to
\begin{equation}\label{equit}
\begin{split}
x_0^2+\dots+x_3^2=N,
\\
x_i\equiv a_i \mod 2q,
\end{split}
\end{equation}
where $a_2=a_3=0$. If there is no solution to \eqref{equit}, then it returns ``No solution''.
\end{thm}
 By Theorem~\ref{reductionthm}, the algorithm in Theorem~\ref{diagliftt} gives the navigation algorithm described   in Theorem~\ref{mainthm}.  
\begin{rem}\label{Npprem}
It is possible to generalize our  polynomial-time algorithm for any  $d\geq2$, by assuming a variant of $(*)$ and two coordinates of the congruence conditions are  zero. On the other hand,  by assuming GRH and Cramer conjecture, Corollary~\ref{cooor}  implies  that the complexity of the optimal strong approximation for a generic point on the  sphere is NP-complete.
%
Hence, by assuming these widely believed arithmetic assumptions, Corollary~\ref{cooor} essentially implies that finding the shortest possible path between a generic pair of vertices in LPS Ramanujan graphs is NP-complete.
\end{rem}
%
 
%
%
%

%
%

 \subsection{Quantitative bounds on the size of the output}\label{quant}
In this section, we give a correspondence between  the  diagonal vertices of  $X^{p,q}$ and the index $q$ sublattices of $\mathbb{Z}^2$. Next, we relate the graph distance between the diagonal vertices (that is a diophantine exponent by Theorem~\ref{reductionthm} ) to the length of the shortest vector of the corresponded sublattice.

 Let $v\in \begin{bmatrix}a+ib & 0\\0 & a-ib \end{bmatrix} \in X^{p,q} $ be a diagonal vertex, and  let $L_v$ be the sublattice of  $\mathbb{Z}^2$ defined by the following congruence equation: 
$$ax+by\equiv 0  \text{ mod } q.$$
Let $\{u_1,u_2\}$ be the Gauss reduced basis for $L_v$, where $u_1$ is a shortest vector in $L_v$. In the following theorem, we relate the graph distance of $v$ from the identity to the norm of  $u_1$. 
\begin{thm}\label{correspond}
Assume Conjecture~\ref{cc}. Let $v$, $L_v$ and $\{u_1,u_2\}$ be as above. Suppose that $\frac{|u_2|}{|u_1|}\geq C_{\gamma} \log(2q)^{\gamma} ,$ then the distance of $v$ from the identity is less than 
\begin{equation}\label{bigholes}\lceil 4\log_p(q)-2\log_{p}|u_1| +\log_p(89)\rceil.\end{equation}
Otherwise, the distance of $v$ from the identity vertex is less than
\begin{equation}\label{almost}\lceil 3\log_p(q)+\gamma\log_{p}\log(q)+\log_{p}(C_{\gamma}) +\log_p(89)\rceil.\end{equation}
\end{thm}
\begin{rem}
In Section~\ref{numericss}, we numerically check  that the inequality~\eqref{bigholes} is sharp. 
 In particular,  the diameter of LPS Ramanujan graphs is asymptotically the longest  distance between the diagonal vertices.  Moreover, the above theorem implies~\eqref{almostdiam} and \eqref{diamb} in Theorem~\ref{mainthm}.
We also use this theorem in our algorithm in Theorem~\ref{decomdiag}, in order to avoid the diagonal vertices with long distance from the identity. 
  \end{rem}

\subsection{Further motivations and techniques}

Rabin and  Shallit \cite{Rabin} developed a randomized  polynomial-time algorithm that represents any integer as a sum of four squares.  
The question of representing a prime as a sum of two squares in  polynomial-time has been discussed in \cite{Schoof} and  \cite{Rabin}. Schoof   developed a deterministic  polynomial-time algorithm that represents a prime $p\equiv 1$ mod 4 as a sum of two squares by  $O((\log p)^6)$ operations. We use Schoof's algorithm in our algorithm in Theorem~\ref{diagliftt}.

%
   
    Both Ross-Selinger and our algorithm start with searching for integral lattice points inside a convex region that is defined by a simple system of quadratic inequalities. If the convex region is defined by a system of linear inequalities in a fixed dimension then the general  result of Lenstra \cite{Lenstra} implies this search  is polynomially solvable. We use a variant of Lenstra's argument in the proof of Theorem~\ref{diagliftt}.  An important feature of our algorithm is that it has been implemented,  and it runs and terminates quickly. We give our numerical results in Section \ref{numericss}.

\subsection*{Acknowledgements}\noindent I would like to thank my Ph.D. advisor Peter Sarnak for several insightful and inspiring conversations during the course of this work. Furthermore, I am very grateful for his letter to me which deals with the Archimedean version of Theorem~\ref{nptheorem}. I would like to thank Professor Peter Selinger for providing a public library of his algorithms in~\cite{Selinger}. This material is partially supported by the National Science Foundation under Grant No. DMS-1440140 while the author was in residence at the Mathematical Sciences Research Institute in Berkeley, California, during the Spring 2017 semester.

\section{NP-Completeness}\label{NP} 
In this section, we prove Theorem~\ref{nptheorem} and Corollary~\ref{cooor}. We reduce them to the sub-sum problem, which is well-known to be NP-complete. We begin by stating  the sub-sum problem, and proving some auxiliary  lemmas. The proof of~Theorem~\ref{nptheorem} and Corollary~\ref{cooor} appear at the end of this section.

 Let $t_1, t_2, \dots, t_k$ and $t\in \mathbb{N}$ with $\log(t)$ and $\log(t_i)$ at most $k^{A}.$ 
 
%
\begin{q} Are there $\epsilon_{i}\in \{0,1 \}$ such that 
\begin{equation}\label{subsum}\sum_{j=1}^{k} \epsilon_j t_j=t?\end{equation}
\end{q}
\begin{lem}
By Cramer's conjecture, there exists a polynomial-time algorithm in $k$ that returns a prime number $q\equiv 3$ mod $4$ such\begin{equation}\label{boundd}
 q> 2 k \max_{1 \leq i \leq k}(t_i,t).
 \end{equation}  Alternatively, this task can be done unconditionally by a probabilistic  polynomial-time algorithm  in $k$.
\end{lem}
\begin{proof}
Let $X:=4  k \max_{1 \leq i \leq k}(t_i,t)+3.$
  We  find $q$ by running the primality test algorithm of Agrawal, Kayal and  Saxena  \cite{primaty}  on the arithmetic progression   $X, X+4, \dots$. By Cramer's conjecture this search terminates in $O(\log(X)^2)$ operations. 
  
  Alternatively,  we  pick a random number between $[X,2X]$ and check by a the primality test algorithm if the number is prime. The expected time of the operations is $O(\log(X))$.
   \end{proof}
  
   Let
 \begin{equation}\label{ss}
 s:=(q-1)t +\sum_{i=1}^{ k} t_i.
 \end{equation}
 By a simple change of variables, solving equation~\eqref{subsum} is equivalent to solving 
 \begin{equation}\label{frob}
 \sum_{j=1}^{k}\xi_{j}t_{j}=s,
 \end{equation}
where $\xi_{j}\in \{1,q   \}$. 

 Let $\mathbb{F}_{q^2}$ be the finite field with $q^2$ elements. 
 \begin{lem}
  By assuming GRH, there exists  a deterministic  polynomial-time algorithm in $\log q$  that returns a finite subset $H\subset \mathbb{F}_{q^2}^*$ of size $O((\log q)^{8+\epsilon})$ such that $H$ contains at least  a generator for the cyclic  multiplicative group  $\mathbb{F}_{q^2}^*$.  Alternatively, this task can be done unconditionally by a probabilistic polynomial algorithm.
 \end{lem}
 \begin{proof}
Since $q\equiv 3$ mod 4, $\mathbb{Z}/q\mathbb{Z}[i]$ is isomorphic to $\mathbb{F}_{q^2}$, where $i^2=-1$. 
  By Shoup's result  \cite[Theorem 1.2]{Shoup},  there is a primitive roots of unity $g=a+bi \in \mathbb{Z}/q\mathbb{Z}[i]$ for the finite field with $q^2$ elements such that $a$ and $b$ has an integral lift of size $O(\log(q)^{4+\epsilon} ) $ for any $\epsilon >0$ result. Hence, the reduction of $H:\{a+bi: |a|,|b|\leq \log(q)^{4+\epsilon} \}$ mod $q$ has the desired property.    
  
  Alternatively, this task can be done unconditionally by a probabilistic polynomial algorithm. Because, the density of primitive roots of unity in $\mathbb{F}_{q^2}^*$ is $\varphi (q^2-1)/(q^{2}-1)$, where $\varphi$ is   Euler's totient  function.  The ratio $\varphi (q^2-1)/(q^{2}-1)$  is well-known to be $O(\log \log q )$. 
%
%
\end{proof}

Next, we take an element $g\in H$ (not necessarily a generator), and  let 
\begin{equation}\label{a,b}
\begin{split}
a_j+ib_j:=g^{t_j} \text{  for  }  1 \leq j \leq k,
\\
a+ib:=g^{s},
\end{split}
\end{equation}
where $t_j$ are given in the subsum problem~\eqref{subsum}, $s$ is defined by equation~\eqref{ss} and $a_j, b_j,a,b \in \mathbb{Z}/q\mathbb{Z}$. Next, we find gaussian primes $\pi_j\in \mathbb{Z}[i]$ such that 
\begin{equation}\label{ppp}\pi_j\equiv a_j+i b_j \text{  mod }  q. \end{equation}
Again this is possible  deterministically by Cramer's conjecture. Alternatively, we choose a random integral point $(h_1,h_2)\in [X,2X]\times  [X,2X]$,  and check by a  polynomial-time primality test algorithm if $(h_1q+a_j)^2+(h_2q+b_j)^2$ is prime in $\mathbb{Z}$. Set $p_i:=|\pi_j|^2$ that is a prime in $\mathbb{Z}$ and define

\begin{equation}\label{N}
N:=\prod_{j=1}^k p_j.
\end{equation}
Consider the following diophantine equation 
\begin{equation}\label{diop}
\begin{split}
&X^2+Y^2=N,
\\
&X\equiv a \text{ and } Y\equiv b \text{ mod } q,
\end{split}
\end{equation}
where $a,b$ are defined in equation~\eqref{a,b} and $N$ in \eqref{N}. Our theorem is a consequence of the following lemma.
\begin{lem}\label{primitive}
Assume that $g\in  \mathbb{F}_{q^2}^*$ is a generator. An integral solution $(X,Y)$ to the diophantine equation ~\eqref{diop} gives a solution $(\xi_1,\dots,\xi_k)$ to the equation~\eqref{frob} in  polynomial-time in $\log(q)$.
\end{lem}
\begin{proof}[Proof of Lemma~\ref{primitive}]  Assume that the equation~\eqref{diop} has an integral solution $(a_0,a_1)$. $A+Bi$ factors uniquely  in $\mathbb{Z}[i]$, and we have
\begin{equation}\label{fact}A+iB=\pm i\prod_{j=1}^{k}\pi_{j}^{\epsilon_j},\end{equation}
where $\epsilon_{j}\in\{0,1 \}$, and $\pi_{j}^{0}=\pi_{j}$ and  $\pi_{j}^{1}=\bar{\pi}$ (the complex conjugate of $\pi_j$). We consider the above equation mod $q$. Then
$$A+iB\equiv \pm i\prod_{j=1}^{k}\pi_{j}^{\epsilon_j} \text{  mod } q.$$
By the congruence condition~\eqref{diop}, $A+iB\equiv a+ib$  mod   $q$, and  by~\eqref{ppp}, $\pi_{j}^{0}\equiv a_j+i b_j $ and  $\pi_{j}^{1}\equiv a_j-i b_j $  mod   $q.$
By~\eqref{a,b}, we obtain 
\begin{equation*}
g^{s}\equiv \prod_{j=1}^{k} g^{\xi_{j} t_{j}},
\end{equation*}
where $\xi_{j}=1$ if $\epsilon_{j}=0$ and $\xi_{j}=q$ if $\epsilon_{j}=1$. Therefore we obtain the following congruence equation
$$ \sum_{j=1}^k \xi_{j} t_{j} \equiv s \text{  mod  }  q(q-1). $$
By the inequality~\eqref{boundd} and the definition of $s$ in~equation~\eqref{ss}, we deduce that 
$$\sum_{j=1}^k \xi_{j} t_{j} =s.$$
This completes the proof of our lemma.\end{proof}
\begin{proof}[Proof of Theorem~\ref{nptheorem}]
Our theorem is a consequence of Lemma~\ref{primitive}. For every $g\in H$ we apply Lemma~\ref{primitive} and check if $(\xi_1,\dots,\xi_k)$ is a solution to~\eqref{frob}. Since the size of $H$ is $O(\log(q)^{8+\epsilon})$ and it  contains at least a primitive roots of unity, we  find a solution to the equation~\eqref{frob} in  polynomial-time. This concludes our theorem.

\end{proof}
\noindent Finally,  we give a proof for Corollary~\ref{cooor}. 
\begin{proof}[Proof of Corollary~\ref{cooor}] We prove this corollary by induction on $d$. The base case $d=2$ follows from Theorem~\ref{nptheorem}.   It suffices to reduce the task with $d$ variables to a similar task with $d+1$ variables in  polynomial-time. The task with $d$ variables is to accept $(N,q,a_1, \dots,a_d)$ such that the following diophantine equation has a solution  

 \begin{equation}\label{deqq}
\begin{split}
&X_1^2+\dots +X_d^2=N,
\\
&X_1\equiv a_1 \dots X_d\equiv a_d \text{ mod } q.
\end{split}
\end{equation}
%
%
We proceed by  taking auxiliary parameters   $0 \leq t, m \in \mathbb{Z}$ such that $N<q^{2t}$,   $m\leq (1/3) q^{2t+1}$ and $\gcd(m,q)=1$. We  consider the following diophantine equation
\begin{equation}\label{3eq}
\begin{split}
&X_1^2+ \dots X_d^2+X_{d+1}^2= m^2+q^{2t}N,
\\
&X_1\equiv q^ta_1, \dots  ,X_d\equiv a_dq^t  \text { and }  X_{d+1}\equiv m\text{ mod } q^{t+1}.
\end{split}
\end{equation}
Assume that $(X_1, \dots, X_{d+1})$ is a solution to the above equation. Then
$$X_{d+1}\equiv \pm m \text{  mod }  q^{2t+1}.$$
 Since $m\leq (1/3) q^{2t+1}$, either $X_{d+1}=\pm m \text{ or }|X_{d+1}|\geq (2/3)q^{2t+1}$. If $|X_{d+1}|\geq (2/3)q^{2t+1}$,  since $m\leq (1/3) q^{2t+1}$ and $N<q^{2t}$, 
 $$X_{d+1}^2> m^2+q^{2t}N. $$
This contradicts with equation~({\ref{3eq}}). This shows that $X_{d+1}=\pm m$. Hence,  the integral solutions to the diophantine equation~\eqref{3eq} are of the form $(q^t X_1, \dots, q^t X_d, \pm m)$ such that $(X_1, \dots, X_d)$ is a solution to the equation~\eqref{deqq}. By our induction assumption, this problem is NP-complete, and we conclude our corollary.

\end{proof}

\section{Algorithm}\label{alg}
\subsection{Proof of Theorem~\ref{diagliftt}}
  In this section, we prove Theorem~\ref{diagliftt}, which is the main ingredient   in   the navigation algorithms in  Theorem~\ref{mainthm} and Theorem~\ref{decomdiag}.

 Let $(x_0,\dots,x_3)\in\mathbb{Z}^4$ be a solution  to the equation~\eqref{equit}.
We change the variables to $(t_0,\dots,t_3)\in\mathbb{Z}^4$, where $x_i=2t_iq+a_i,$
and $|a_i| \leq q$. Hence,
\begin{equation}\label{main}\frac{N}{4q^2}-(t_0+a_0/2q)^2-(t_1+a_1/2q)^2= t^2_2 +t^2_3.\end{equation}
Let $Q(t_0,t_1):=\frac{N}{4q^2}-(t_0+a_0/2q)^2-(t_1+a_1/2q)^2.$ Recall the definition of $A_{Q,r}$ from~\eqref{Afr}, where $r>0$ is some real number. By conjecture~\ref{cc}, if $ |A_{Q,r}|> C_{\gamma}(\log N)^{\gamma}$ then the equation~\eqref{main} has a solution, where $(t_0,t_1)\in A_{Q,r}.$ 

First, we give a parametrization of $(t_0,t_1)\in \mathbb{Z}^2$, where $Q(t_0,t_1)\in \mathbb{Z}.$
Let $k:=\frac{N-a_0^2-a_1^2}{4q}$.  Since $a_0^2 +a_1^2 \equiv N  \text{  mod } 4q$, $k\in\mathbb{Z}$. 
%
By~\eqref{main},
\begin{equation}\label{newconj} a_0t_0+a_1t_1\equiv k  \text{ mod } q.\end{equation}
 Without loss of generality, we assume that  $a_0\neq 0$ mod $q$.  Then $a_0$ has an inverse mod $q$, and $(ka_0^{-1},0)$ is a solution for the congruence equation~\eqref{newconj}. We lift  $(ka_0^{-1},0)\in\big(\mathbb{Z}/q\mathbb{Z}\big)^2$ to the integral vector $(c,0)\in \mathbb{Z}^2$ such that 
$$c\equiv ka_0^{-1} \text{ mod } q \text{ and } |c|<(q-1)/2. $$
The integral solutions of equation~\eqref{newconj} are the translation of the integral solutions of the following  homogenous equation by vector $(c,0)\in \mathbb{Z}^2$
\begin{equation}\label{conj2} a_0t_1+a_1t_1\equiv 0  \text{ mod } q.\end{equation}
The integral solutions to equation~({\ref{conj2}}) form a lattice of co-volume $q$ that is spanned by the integral basis $\{v_1,v_2  \}$ where
$$v_1:=(q,0), \text{ and }  v_2:=(-a_1a_0^{-1},1).$$
We apply  Gauss reduction algorithm on the basis $\{v_1,v_2\}$  in order to find an almost orthogonal basis $\{u_1,u_2\}$ such that  
\begin{equation}\label{ortho}
\begin{split}
\text{span}_{\mathbb{Z}}\langle v_1,v_2 \rangle =\text{span}_{\mathbb{Z}}\langle u_1, u_2 \rangle,
\\
|u_1|<|u_2|,
\\
\langle u_1,u_2\rangle \leq  (1/2)\langle u_1,u_1\rangle. 
\end{split}
\end{equation}
where  $\text{span}_{\mathbb{Z}}\langle v_1,v_2 \rangle:=\{xv_1+yv_2: x,y \in \mathbb{Z}     \}$ and $\langle u_1,u_2\rangle\in \mathbb{R}$ is the dot product of $u_1$ and $u_2$. 
%
%
 Let $u_0$ be a shortest integral vector that satisfies the equation~\eqref{newconj}. We write $(c,0)$ as a linear combination of $u_1$ and $u_2$ with coefficients in $(1/q) \mathbb{Z}$
$$(0,c)=(h_1+r_1/q)u_1+ (h_2+r_2/q)u_2,$$
where $0 \leq r_1,r_2\leq q-1$.  Note that $u_0$ is one of the following 4 vectors
$$ \big(r_1/q-\{0,1 \}   \big) u_1 + \big(r_2/q-\{0,1 \}  \big) u_2.$$ 
By triangle inequality,  $$|u_0|<|u_2|.$$ We parametrize the integral solutions $(t_0,t_1)$ of~\eqref{newconj} by:
\begin{equation}\label{param}(t_0,t_1)=u_0+xu_1+yu_2,\end{equation}
where $x,y \in \mathbb{Z}.$ Let $u_0=(u_{0,0},u_{0,1})$, $u_1=(u_{1,0},u_{1,1})$ and $u_2=(u_{2,0},u_{2,1})$. Since $u_1$ and $u_2$ are solutions to~\eqref{conj2} and $u_0$ is a solution to~\eqref{newconj}, 
\begin{equation*}
\begin{split}
u_0^{\prime}:= \frac{k-a_0u_{0,0}-a_1u_{0,1}}{q} \in\mathbb{Z},
\\
u_1^{\prime}:=\frac{a_0u_{1,0}+a_1u_{1,1}}{q}\in\mathbb{Z},
\\
u_2^{\prime}:=\frac{a_0u_{2,0}+a_1u_{2,1}}{q}\in\mathbb{Z}.
\end{split}
\end{equation*}
 Let 
 \begin{equation}\label{F(x)}F(x,y):=u_0^{\prime}-xu_1^{\prime}-yu_2^{\prime}-(u_{0,1}+xu_{1,1}+yu_{2,1})^2 -(u_{0,2}+xu_{1,2}+yu_{2,2})^2. \end{equation}
By \eqref{param},
\begin{equation*}
F(x,y)=Q(t_0,t_1).
 \end{equation*}
 Hence, $Q(t_0,t_1)\in \mathbb{Z}$ for $(t_0,t_1)\in \mathbb{Z}^2,$ if and only if $(t_0,t_1)=u_0+xu_1+yu_2$ for some $(x,y)\in\mathbb{Z}^2.$

Next, we list  all the integral points $(x,y)$ such that $F(x,y)$ is positive. 
\begin{lem}\label{boxlem}
Assume that  $\frac{\sqrt{N}}{q|u_2|}\geq 14/3$. Let $F(x,y)$ be as above. Let $A:=\sqrt{N}/(2q|u_1|)-1$ , $B:=\sqrt{N}/(2q|u_2|)-1$ and 
\begin{equation}\label{box} C:= [-A,A]\times[-B,B]. \end{equation} Then  $F(x,y)$ is positive for every $(x,y)\in C$ and negative outside $10\times C$.
\end{lem}
 \begin{proof}
  Recall that $(t_0,t_1)=u_0+xu_1+yu_2$ and
$$ F(x,y)=  N/4q^2-(t_1+a_0/4q)^2 - (t_2+a_1/2q)^2,$$
where $|a_0/2q|<1/2 \text{ and } |a_1/2q| < 1/2$. Hence,  if  $|(t_0,t_1)|<(\sqrt{N}/q) -1 $, then $F(x,y)>0$, and if $|(t_0,t_1)|>(\sqrt{N}/q) +1,$ then $F(x,y)<0$. By the triangle inequality 
\begin{equation*}
\begin{split}
|(t_0,t_1)|=|u_0+xu_1+yu_2|
\leq |u_0|+|x||u_1|+|y||u_2|.
\end{split}
\end{equation*}
Since $|u_0|<|u_2|$ then $|(t_0,t_1)|\leq |x| |u_1|+(1+|y|)|u_2|.$ Let $A$ , $B$ and $C$ be as in~\eqref{box}. 
Then for every $(x,y)\in [-A,A]\times[-B,B]$, we have 
$$
|x| |u_1|+(1+|y|)|u_2|\leq (\sqrt{N}/q) -|u_1| < (\sqrt{N}/q) -1.
$$
Hence, $F(x,y)>0$ if $(x,y)\in [-A,A]\times[-B,B]$. Next, we show that $F$ is negative outside $10\times C.$ By almost orthogonality conditions~\eqref{ortho}, we obtain the following lower bound 
\begin{equation}\label{lowerbd}(|x|/2)|u_1|+(|y|/2-1)|u_2| \leq |u_0+xu_1+yu_2|.\end{equation}
The above  inequality implies that if  $x \geq 10A$, then 
 $$|(t_0,t_1)|=|u_0+xu_1+yu_2| > \sqrt{N}/q+ \big(3\sqrt{N}/q|u_1|-10\big)|u_1|/2.$$ 
We assume that $\frac{\sqrt{N}}{q|u_2|}\geq 14/3$ and $1<|u_1|<|u_2|$, then $|(t_0,t_1)|> \sqrt{N}/q+1 $ and hence $F(x,y)$ is negative. Similarly, if  $y \geq 10 B $ then
  $$|(t_0,t_1)|=|u_0+xu_1+yu_2| > \sqrt{N}/q+ \big(3\sqrt{N}/(2q|u_2|)-6\big)|u_2|.$$ 
 Since $\frac{\sqrt{N}}{q|u_2|}\geq 14/3$, it follows that  $|(t_0,t_1)|> \sqrt{N}/q+1.$ Hence,  $F(x,y)$ is negative. Therefore, if $(x,y)\notin 10 \times  C$, then $F(x,y)$ is negative.  This concludes our lemma. 
\end{proof}
In the following lemma, we consider the remaining case, where $\frac{\sqrt{N}}{q|u_2|}\leq 14/3.$ 

\begin{lem}\label{line}
Assume that  $\frac{\sqrt{N}}{q|u_2|}\leq 14/3$ and $F(x,y)>0$ then $|y|\leq 13.$ 
\end{lem} 
\begin{proof} Since  $F(x,y)>0$ from the first line of the proof of Lemma~\ref{boxlem}, it follows that $|(t_0,t_1)|< (\sqrt{N}/q) +1.$  From the the inequality~\eqref{lowerbd}, we have  
$$
(|y|/2-1)|u_2| \leq |(t_0,t_1)| \leq \sqrt{N}/q +1.
$$
Hence, 
$$
|y|\leq  \frac{2\sqrt{N}}{q|u_2|} +4 < 14.$$
Since $y$ is an integer, we conclude the lemma. 
\end{proof}

\begin{proof}[Proof of Theorem~\ref{diagliftt}]
 Assume that  $\frac{\sqrt{N}}{q|u_2|}\geq 14/3.$  By Lemma~\ref{boxlem},  $F(x,y)$ is positive inside box $C$ that is defined in~\eqref{box}. We list  $(x,y)\in C$ in the order of their distance from the origin. If possible  we represents $F(x,y)$  as a sum of two squares by the following polynomial-time  algorithm. We factor $F(x,y)$ into primes by the polynomial algorithm for factoring integers in $(*)$. Next, by Schoof's algorithm \cite{Schoof}), we write every prime number as a sum of two squares. If we succeed, then we find an integral solution to the equation~\eqref{equit}, and this concludes the theorem. 
 
  
  If the size of box $C$ that is $A\times B > C_{\gamma}\log(N)^{\gamma},$ then by Conjectrue~\ref{cc} we find a pair $(x,y)$ such that $F(x,y)$ is a sum of two squares in less than $ O(\log(q)^{O(1)})$ steps, and the above algorithm terminates. Otherwise,  $A\times B < C_{\gamma}\log(M)^{\gamma}$.  By Lemma~\ref{boxlem}, $F(x,y)$ is negative outside box $10C$ and since the size of this box is $O(\log(q)^{\gamma})$ we check all points inside box $10C$ in order to represent $F(x,y)$ as a sum of two squares. If we succeed to represent $F(x,y)$ as a sum of two squares then we find an integral solution to equation~\eqref{equit}. Otherwise, the equation~\eqref{equit} does not have any integral solution. This concludes our theorem if $\frac{\sqrt{N}}{q|u_2|}\geq 14/3$. 
  
  Finally, assume that $\frac{\sqrt{N}}{q|u_2|}\leq 14/3$ then by Lemma~\ref{line}, we have $|y| \leq 13.$ We fix $y=l$ for some $|l|<13.$ We note that by equation~\eqref{F(x)}, 
$$
F(x,l)=A x^2+B x+ C
$$
for some $A, B,C\in \mathbb{Z} $. We list  $x\in \mathbb{Z}$ such that $F(x,l)>0$ and then proceed similarly as in the forth line of the first paragraph of the proof. This concludes our theorem.  
%
%
%
%
%
%
%

\end{proof}

\subsection{Distance of diagonal vertices from the identity}
  In this section, we give a proof of Theorem~\ref{correspond}. Then, we give bounds on the size of the outputs in Theorem~\ref{mainthm} and Theorem~\ref{decomdiag}. Recall the notations while formulating  Theorem~\ref{correspond}. 

\begin{proof}[Proof of Theorem~\ref{correspond}]
%
 
 We proceed by proving \eqref{bigholes}. Assume that $$|u_2|\geq C_{\gamma} \log(q)^{\gamma} |u_1|.$$  Let 
 \begin{equation}\label{hbound}
 h:=\lceil 4\log_p(q)-2\log_{p}|u_1| +\log_{p}(89) \rceil.
 \end{equation}
 We show that there exists a path from $v$ to the identity of length $h$.  By our assumption $p$ is a quadratic residue mod $q$. We denote the square root of $p$ mod $q$ by $\sqrt{p}.$  Set 
 \begin{equation*}
 \begin{split}
 A:=a\sqrt{p}^h \text{  mod }  2q,
 \\
 B:=b\sqrt{p}^h \text{  mod  } 2q.
 \end{split}
 \end{equation*}
By Theorem~\ref{reductionthm}, there exists a path of length $h$ from $v$ to the identity if and only the following diophantine equation has an integral solution $(t_1,t_2,t_3,t_4)$
 \begin{equation}\label{newdiag}(2t_1q+A)^2 +  (2t_2q+Bt)^2 + (2t_3q)^2 +(2t_4q)^2=p^{h}. \end{equation}
In Theorem~\ref{diagliftt}, we developed a  polynomial-time algorithm for finding its integral solutions $(t_1,t_2,t_3,t_4)$. 
We defined the associated binary quadratic form $F(x,y)$ as defined in equation~\eqref{F(x)}. By Lemma~\ref{boxlem},  $F(x,y)$ is positive inside the box $[-A,A]\times[-B,B]$ where $A:=\sqrt{p^h}/(4q|u_1|)$ , $B:=\sqrt{p^h}/(4q|u_2|)-1$. By the definition of $h$ in equation~\eqref{hbound}, we have
\begin{equation}
p^{h}\geq \frac{89q^4}{|u_1|^2}.
\end{equation}
By the above inequality 
\begin{equation}\label{Bin}
B\geq  \frac{\sqrt{89}q^2}{4q|u_1||u_2|}-1.
\end{equation}
Since $\{u_1,u_2\}$ is an almost orthogonal basis for a co-volume $q$ lattice then the angle between $u_1$ and $u_2$ is between $\pi/3$ and $2\pi/3$. Hence,
\begin{equation}\label{area}|u_1||u_2|\leq 2q/\sqrt{3}.\end{equation}
We use the above bound on $|u_1||u_2|$ in inequality~\eqref{Bin}, and derive
$$B\geq  \frac{\sqrt{3*89}}{8}-1 >1.$$
Next, we give a lower bound on $A$. Note that 
$$A \geq \frac{|u_2|}{|u_1|}B.$$
By our assumption $\frac{|u_2|}{|u_1|}\geq C_{\gamma} \log(2q)^{\gamma}$, hence
$$A \geq C_{\gamma} \log(2q)^{\gamma} B.$$
Since $B>1$, $$AB\geq C_{\gamma} \log(2q)^{\gamma}.  $$
By Conjecture~\ref{cc} and Theorem~\ref{diagliftt}, our algorithm returns an integral solution $(t_1,t_2,t_3,t_4)$ which gives rise to a path of length $h$ from $v$ to the identity. This concludes the first part of our theorem.

 Next, we assume that $|u_1|\leq |u_2|\leq C_{\gamma} \log(2q)^{\gamma} |u_1|$. Let 
 \begin{equation}\label{kdef}h^{\prime}:=\lceil 3\log_p(q)+\gamma\log_{p}\log(q)+\log_{p}(C_{\gamma}) +\log_p(89)\rceil.\end{equation}
We follow the same analysis as in the first part of the theorem. First, we  give a lower bound on $B:=\sqrt{p^{h^{\prime}}}/(4q|u_2|)-1$. By the definition of $h^{\prime }$ in equation~\eqref{kdef}, we derive
\begin{equation}\label{pkin}
p^{h^{\prime}}\geq  89C_{\gamma}\log(q)^{\gamma}q^3.
\end{equation}
 We multiply both sides of  $|u_2|\leq C_{\gamma} \log(q)^{\gamma} |u_1|$ by $|u_2|$ and use the inequality~\eqref{area} to obtain
$$|u_2|^2 \leq C_{\gamma} \log(q)^{\gamma}2q/\sqrt{3}. $$ 
 By the above inequality, definition of $B$ and inequality~\eqref{pkin}, we have
 \begin{equation*}B=\sqrt{p^{h^{\prime}}}/(4q|u_2|)-1\geq \frac{\sqrt{89}}{4\sqrt{2/\sqrt{3}}}-1\geq 1.\end{equation*}
 Hence, $$B\geq \sqrt{p^{h^{\prime}}}/(8q|u_2|).$$ Next, we use the above inequality and inequality~\eqref{area} and \eqref{pkin} to give a lower bound on $AB$.
 \begin{equation}
 \begin{split}
 AB&\geq \frac{p^{h^{\prime}}}{32q^2|u_1||u_2|} 
 \\
 &\geq \frac{ 89\sqrt{3}C_{\gamma}\log(q)^{\gamma}q^3}{64q^3}
 \\
 &> C_{\gamma}\log(q)^{\gamma}.
 \end{split}
 \end{equation}
 By Conjecture~\ref{cc} and  Theorem~\ref{diagliftt}, our algorithm returns an integral solution $(t_1,t_2,t_3,t_4)$ which gives rise to a path of length $h^{\prime}$ from $v$ to the identity. This concludes our theorem. 
\end{proof}
%

Finally, we prove \eqref{almostdiam} in Theorem~\ref{almostdiam}. We briefly, explain the main idea.  We normalize the associated co-volume $q$ lattices $L_v$, so that they have co-volume 1. These normalized lattices are parametrized by points in  $\text{SL}_2(\mathbb{Z})\backslash \mathbb{H}$, and  it is well-known that they are equidistributed in $\text{SL}_2(\mathbb{Z})\backslash \mathbb{H}$ with respect to the hyperbolic measure $\frac{1}{y^2} (dx^2+dy^2).$
%
%
It follows from this equidistribution and Theorem~\ref{correspond} that  the distance of a typical diagonal matrix from the identity vertex is $\log(|X^{p,q}|)+O(\log \log (|X^{p,q}|))$. The diagonal points with distance $4/3\log(|X^{p,q}|)+O(\log \log (|X^{p,q}|))$ from the identity are associated to points $x+iy\in \mathbb{H}$ with $y$  as big as $q$. 
%

\begin{proof}[Proof of  \eqref{almostdiam} in Theorem~\ref{almostdiam}] Let $v$ be a diagonal vertex with distance $h$ from the identity vertex, where $$ h\geq  \lceil 3\log_p(q)+\gamma\log_{p}\log(q)+\log_{p}(C_{\gamma}) +\log_p(89) \rceil. $$  Let $L_v$ be the associated lattice of co-volume 
$q$ and $\{u_1,u_2 \}$ be an almost orthogonal basis for $L_v$.  By Theorem~\ref{correspond}, the distance of $v$ from the identity is less than
$$\lceil 4\log_p(q)-2\log_{p}|u_1| +\log_p(89)\rceil.$$ 
Therefore,
$$ h \leq \lceil 4\log_p(q)-2\log_{p}|u_1| +\log_p(89)\rceil.$$
Hence, 
\begin{equation}\label{u1}|u_1|^2 \leq 89q^4/p^{(h-1)}.\end{equation}
Next, we count the number of lattices of co-volume $q$ inside $\mathbb{Z}^2$ such that the length of the shortest vector is smaller than  $r\leq (1/2) \sqrt{q}$. Let $L\subset \mathbb{Z}^2$ be a lattice of co-volume $q$ such that $L$ contains a vector of length smaller than $(1/2) \sqrt{q}$. It is easy to check that $L$ contain unique vectors $\pm v:=\pm(a_0,a_1)$ such they have the shortest length among all vectors inside $L$. Since $q$ is prime this vector is primitive i.e. $\gcd(a_0,a_1)=1$. On the other hand, the lattice is uniquely determined by $\pm v:=\pm(a_0,a_1)$, namely $L$ is the set of all integral points $(x,y)\in \mathbb{Z}^2$ such that 
$$ax+by\equiv 0 \text{ mod } q.$$
Therefore, the problem of counting the lattices of co-volume $q$ with shortest vector smaller than $r$ is reduced to counting the projective primitive integral vectors of length smaller than $r$. The main term of this counting is 
\begin{equation}\label{count}1/2 \zeta(2)^{-1} \pi r^2= \frac{3}{\pi}r^2.\end{equation}
By inequality~\eqref{u1} and~\eqref{count}, we deduce that the number of diagonal vertices  
with graph distance at least $h$ from the identity in LPS Ramanujan graph $X^{p,q}$ is less than  
$$89q^4/p^{(h-1)}.$$
This concludes Theorem~\ref{mainthm}.

\end{proof}

\subsection{Algorithm for the diagonal decomposition}\label{diagsection}
%
\begin{proof}[Proof of \eqref{deter} in Theorem~\ref{decomdiag}]  Let $M:=\begin{bmatrix}a & b\\ c &d   \end{bmatrix}\in  {\rm PSL}_2(\mathbb{Z}/q\mathbb{Z})$ be any element.   By \cite[Lemma~1]{Petit2008}, there exists a polynomial-time algorithm that expresses $M$ as:
$$
M=D_1s_1D_2s_2D_3s_3D_4,
$$ 
where $D_i$ are diagonal matrices for $1\leq i\leq 4$ and $s_j$ are LPS generators for $1\leq j\leq 3.$ By Theorem~\ref{mainthm} and assuming $(*)$, we write each $D_i$  as a product of at most  $4/3\log_{p}|X^{p,q}|+O(1)$ LPS generators  in polynomial time.  Therefore, we find a path of size at most  $\frac{16}{3}\log_{k-1}|X^{p,q}|+O(1)$ from the identity to $M.$ This concludes \eqref{deter}.
\end{proof}

Let 
$$
D:=\Big\{\begin{bmatrix}a & 0\\0 &a^{-1}   \end{bmatrix}\in  {\rm PSL}_2(\mathbb{Z}/q\mathbb{Z})   \Big\}, \text{ and }R:=\Big\{\begin{bmatrix}a & b\\ b &a   \end{bmatrix}\in  {\rm PSL}_2(\mathbb{Z}/q\mathbb{Z})   \Big\}. $$
%
Define $d_{\alpha}:=\begin{bmatrix}\alpha & 0 \\ 0&\alpha^{-1}  \end{bmatrix},$ and $ r_{a,b}:= \begin{bmatrix}a & b \\ a&b  \end{bmatrix}.$ 
By the correspondence~\eqref{cores} between ${\rm PSL}_2(\mathbb{Z}/q\mathbb{Z})$ and the units of $\mathbb{H}(\mathbb{Z}/q\mathbb{Z}) $, $D$ and $R$ are associated to:
$$
\tilde{D}:=\{a+bi:a,b\in \mathbb{Z}/q\mathbb{Z}, a^2+b^2=1  \}, \text{ and }\tilde{R}:=\{a+bj:a,b\in \mathbb{Z}/q\mathbb{Z}, a^2+b^2=1  \}.
$$
By Theorem~\ref{diagliftt}, there is a polynomial-time  algorithm that finds the shortest possible path between the identity and vertices in $D$ or  $R.$ Let $D_1\subset D$ and $R_1\subset R$ be the subset of vertices where their distances from the identity is less than  $\log_p|X^{p,q}|+O(\log_p\log(|X^{p,q}|)).$ By \eqref{almostdiam}, in Theorem~\ref{mainthm}, 
$$
|R_1|\geq 99\% |R| \text{ and } |D_1|\geq 99\% |D|.
$$
 Let $Y:=D_1R_1D_1\subset {\rm PSL}_2(\mathbb{Z}/q\mathbb{Z}).$
\begin{lem} We have
$$
|Y|\geq 10\% |{\rm PSL}_2(\mathbb{Z}/q\mathbb{Z})|.
$$
\end{lem}
\begin{proof}
Let $g\in Y.$ Then,
$$
g=d_{\alpha} r_{a,b} d_{\beta}
$$
for some $a,b, \alpha,\beta\in \mathbb{Z}/q\mathbb{Z}.$ We give an upper bound on the number of different ways of expressing $g$ as $d_{\alpha} r_{a,b} d_{\beta}$, where $ab\neq0.$ Suppose that  $d_{\alpha} r_{a,b} d_{\beta}=d_{\alpha^{\prime}} r_{a^{\prime},b^{\prime}} d_{\beta^{\prime}}.$ Then, it follows that $(\alpha^{-1}\alpha^{\prime})^2=(\beta^{-1}\beta^{\prime})^2=\pm1.$ This shows that $g$ has only 2 representations as $d_{\alpha} r_{a,b} d_{\beta}$.  There are only two elements of $R$ with $ab=0$, which are $R_{1,0}$ and $R_{0,i}.$ Since $q$ is a prime,  $|D|=|R|=\frac{q-1}{2}$, and $|{\rm PSL}_2(\mathbb{Z}/q\mathbb{Z})|=\frac{(q-1)q(q+1)}{2}$.  Therefore, 
$$
|Y|\geq  99\%^3\frac{(q-1)(q-5)(q-1)}{16}.
$$
This concludes our lemma. 
\end{proof}
\begin{lem}\label{Ylem}
Let $g\in X^{p,q}.$ There exists a polynomial-time algorithm in $\log q$ that returns a short path of size at most  $3\log_{k-1}|X^{p,q}|+ O(\log \log (|X^{p,q}|))$ form the identity to $g$, if $g\in Y$. Otherwise, it returns ``Not in Y''. 
\end{lem}
\begin{proof}
Let $g=\begin{bmatrix}g_{1,1} & g_{1,2}\\ g_{2,1} &g_{2,2}   \end{bmatrix}.$ First, we check the solubility of $d_{\alpha} r_{a,b} d_{\beta}=g$ for some $\alpha,\beta, a,$ and $b.$ This is equivalent to the following system of equations: 
\begin{equation*}
\begin{bmatrix}
\alpha\beta a & \alpha \beta^{-1}b
\\
 \alpha^{-1} \beta b& \alpha^{-1}\beta^{-1} a
\end{bmatrix}
=\begin{bmatrix}g_{1,1} & g_{1,2}\\ g_{2,1} &g_{2,2}   \end{bmatrix}.
\end{equation*}
It follows that $a^2=g_{1,1}g_{2,2},$  $b^2= g_{1,2}g_{2,1}.$ By the quadratic reciprocity law, we check in polynomial-time algorithm if $g_{1,2}g_{2,1}$ and $g_{1,1}g_{2,2}$ are quadratic residue mod $q$. If either $g_{1,2}g_{2,1}$  or  $g_{1,1}g_{2,2}$ are quadratic non-residue, the algorithm returns ``Not in Y''.   Otherwise,  by the  polynomial-time algorithm for taking square roots in finite fields (e.g  \cite{Adlemann} or \cite{Shanks} ), we find $a$ and $b$. Similarly, we find $\alpha$ and $\beta.$ By Theorem~\ref{mainthm}, we write $d_{\alpha},$ $d_{\beta}$, and $r_{a,b}$ in terms of the LPS generators and check if they are inside $D_1$ and $R_1$ respectively.  This concludes our lemma.

\end{proof}

We cite the following proposition form \cite[Proposition 2.14]{Ellenberg}.
\begin{prop}[Due to Ellenberg,  Michel, Venkatesh]\label{largedev}
 Fix $\epsilon>0$. For any subset $Y\subset X^{p,q}$ with $|Y|>10\%|X^{p,q}|$, the fraction of non-backtracking paths $\gamma$ of
length $2l$ satisfying:
$$
\Big| \frac{|\gamma\cap Y|}{2l+1}-\frac{|Y|}{|G|}   \Big|\geq \epsilon 
$$
is bounded by $c_1 exp(-c_2l)$, where $c_1, c_2$ depend only on $\epsilon.$ 
\end{prop}

\begin{proof}[Proof of \eqref{typical} in Theorem~\ref{decomdiag}]
%
%
It suffices to navigate from the identity to a given vertex $v\in X^{p,q}.$ Recall that $S_p$ is the LPS generator set defined in \eqref{LPSgen}. Let W be the set of all words of length at most $\log\log q$  with letters in $S_p.$ Note that $|W|=O((\log q)^c),$ where $c=\log(p)$ which only depends on the fixed prime $p.$
By Lemma~\ref{Ylem}, if $wv\in Y$  for some $w\in W$, then we find a path that satisfies  \eqref{typical}. By Proposition~\ref{largedev} for $\epsilon=9\%$, it follows that the fraction of the vertices $v$, such that $wv\notin Y$ for every $w\in W$, is less than $c_1\exp(-c_2\log\log q)=O(\log(q)^{-c_1}).$ This concludes our theorem.  

\end{proof}

\section{Numerical results}\label{numericss}

\subsection{Diagonal approximation with V-gates}
In this section, we give some numerical results  on the graph distance between diagonal vertices in $X_{5,q}$ ($V$-gates), which shows that the inequalities~\eqref{almostdiam}~and \eqref{diamb} are sharp.   In particular, we numerically check that the diameter of $X_{5,q}$ is bigger that 
  $(4/3) \log_{5}|X_{5,q}|+O(1).$

 Let $q$ be a prime number and $q\equiv 1,9 \text{ mod } 20$. The LPS generators associated to $p=5$ are called $V$-gates. $V$-gates are the following 6 unitary matrices:
$$V_X^{\pm}:=\frac{1}{\sqrt{5}} \begin{bmatrix} 1 & 2i \\ 2i &1   \end{bmatrix}^{\pm}, \text{    } V_Y^{\pm}:=\frac{1}{\sqrt{5}} \begin{bmatrix} 1 & 2 \\ -2 &1   \end{bmatrix}^{\pm} \text{ and   } V_Z^{\pm}:=\frac{1}{\sqrt{5}} \begin{bmatrix} 1+2i & 0 \\ 0 &1-2i   \end{bmatrix}^{\pm}.$$
Since $q\equiv 1,9 \text{ mod } 20$, then square root of $-1$ and $5$ exist mod $q$ and we denote them by $\sqrt{5}$ and $i$. So, we can realize these matrices inside $PSL_2(\mathbb{Z}/q\mathbb{Z})$.  The Cayley graph of  $PSL_2(\mathbb{Z}/q\mathbb{Z})$ with respect to V-gates is a 6-regular LPS Ramanujan graph. We run our algorithm to find the shortest path in V-gates from identity to a given typical diagonal matrix$\begin{bmatrix}a+bi & 0 \\ 0 & a-bi    \end{bmatrix}\in PSL_2(\mathbb{Z}/q\mathbb{Z}) .$ 
By Theorem~\ref{reductionthm}, a path of length $m$ from identity  to this diagonal element is associated to the integral solution of the following diophantine equation
\begin{equation}\label{cond}
\begin{split}
x^2+y^2+z^2+w^2=5^m
\\
x \equiv \sqrt{5}^m a  \text {  mod  } q, 
\\
y \equiv \sqrt{5}^m b  \text {  mod  } q, 
\\
z\equiv w \equiv 0 \text{ mod } q,
\\
x\equiv 1 \text{ and  } y \equiv z \equiv w \equiv 0 \text{ mod } 2.
\end{split}
\end{equation}
First, our algorithm in Theorem~\ref{diagliftt} finds an integral solution $(x,y,z,w)$ with the least integer $m$ to the equation \eqref{cond}.   Next from the integral solution $(x,y,z,w)$ it constructs a path in the Ramanujan graph by factoring $x+iy+jz+kw$ into $V$-gates. We give an explicit example. Let $q$ be the following prime number with 100 digits:

\begin{equation*}
\begin{split}
q=6513516734600035718300327211250928237178281758494
\\
417357560086828416863929270451437126021949850746381.
\end{split}
\end{equation*}
For the diagonal matrix let $$a=23147807431234971203978401278304192730471291281$$ and 
$$b=1284712970142165365412342134123412341234121234342141234133$$
The first run of our algorithm returns the following  integral lift
$$x+iy+jz+kw$$
where

\begin{equation*}
\begin{split}
x=-3513686405828860927763754940484616687735954403564689113985383253868329887
\\073895129393123529043092607930187858085249975614142765081986624258530038940271
\\
y=3773156548062114482690557548470637380371201820782668326017207890171886678830
\\601870144317232489264867168831689578223312772963262687237828114002146000356
\end{split}
\end{equation*}
\begin{equation*}
\begin{split}z=696150282464006603091186089706225565057448347974579991940267012475009315401865
\\6570861892918415809962375271929963309479306543335375368842987498287311268
\end{split}
\end{equation*}
\begin{equation*}
\begin{split}w=3888519350877870793211628965104035265911619494928178960777970459693109319153422
\\770196318754816019921662119578623310979387405367017752713898473225295568\end{split}
\end{equation*}
and the associated path in the Ramanujan graph by $V$-gates is 

$$Vy  Vz^{-1}  Vx  Vz  Vx  Vx  Vz  Vz  Vx^{-1}  Vx^{-1}  Vz^{-1}  Vx^{-1}  Vz  Vz  Vy  Vz^{-1}  Vz^{-1}  Vz^{-1}  Vy  Vz^{-1}$$
$$Vy^{-1}  Vx  Vx  Vz  Vx  Vy^{-1}  Vx  Vy^{-1}  Vx  Vz^{-1}  Vy  Vx  Vz  Vz  Vx  Vz^{-1}  Vy^{-1}  Vx  Vx  Vz  Vz  Vx  Vx  $$
$$Vz^{-1} Vx  Vx  Vy^{-1}  Vx^{-1}  Vz  Vy  Vx  Vz  Vy  Vx^{-1}  Vy^{-1}  Vy^{-1}  Vz^{-1}  Vy^{-1}  Vz  Vx^{-1}  Vz^{-1}  Vx^{-1}$$
$$Vx^{-1}  Vz  Vy^{-1}  Vx^{-1}  Vz  Vx^{-1}  Vx^{-1}  Vz^{-1}  Vy  Vz  Vz  Vz  Vy  Vz^{-1}  Vx  Vy  Vx^{-1}  Vz^{-1}  Vx^{-1}  Vz^{-1}$$
$$Vx^{-1}  Vx^{-1}  Vz  Vy^{-1}  Vx^{-1}  Vx^{-1}  Vy^{-1}  Vz^{-1}  Vx  Vz^{-1}  Vx^{-1}  Vy  Vy  Vy  Vy  Vy  Vx^{-1}  Vz  Vx^{-1}  Vz  $$
$$Vy  Vx^{-1}  Vx^{-1}  Vy  Vz^{-1}  Vx  Vx  Vz  Vy^{-1}  Vz^{-1}  Vy  Vz  Vx^{-1}  Vx^{-1}  Vy^{-1}  Vz^{-1}  Vy  Vx^{-1}  Vy  Vz^{-1}  Vy$$
$$Vz  Vz  Vx^{-1}  Vx^{-1}  Vy^{-1}  Vx^{-1}  Vz^{-1}  Vx^{-1}  Vy  Vz  Vy  Vy  Vx^{-1}  Vz^{-1}  Vz^{-1}  Vy  Vy  Vx  Vy  Vy  Vz  Vz $$
 $$Vy  Vz Vx  Vz  Vy  Vz  Vx  Vy  Vz^{-1}  Vy  Vx^{-1}  Vz  Vx  Vz^{-1}  Vy^{-1}  Vx  Vx  Vy^{-1}  Vx  Vy  Vx  Vy^{-1}  Vy^{-1}  Vy^{-1} $$
 $$Vz  Vx  Vy^{-1}  Vz  Vx^{-1}  Vz^{-1}  Vx^{-1}  Vx^{-1}  Vz^{-1}  Vz^{-1}  Vy^{-1}  Vx  Vy^{-1}  Vx^{-1}  Vz^{-1}  Vx^{-1}  Vz  Vx  Vz^{-1} $$
 $$Vy^{-1} Vz^{-1}  Vy  Vx  Vz  Vx^{-1}  Vy^{-1}  Vz^{-1}  Vx^{-1}  Vz^{-1}  Vz^{-1}  Vy  Vx^{-1}  Vy^{-1}  Vz^{-1}  Vy  Vz^{-1}  Vx  Vz  Vx$$
 $$Vx  Vy  Vx^{-1}  Vx^{-1}  Vz^{-1}  Vx  Vz  Vy^{-1}  Vz^{-1}  Vz^{-1}  Vy^{-1}  Vy^{-1}  Vy^{-1}  Vx^{-1}  Vx^{-1}  Vy^{-1}  Vz^{-1}$$
 $$Vy  Vx  Vx  Vx  Vy^{-1}  Vx  Vz^{-1}  Vy^{-1}  Vz  Vz  Vy  Vz  Vy  Vz  Vz  Vx  Vx  Vy^{-1}  Vx^{-1}  Vy  Vz^{-1}  Vy^{-1}  Vx^{-1}$$
 $$Vz^{-1}  Vx^{-1}  Vz  Vx  Vy^{-1}  Vx^{-1}  Vx^{-1}  Vy  Vx  Vy  Vx  Vz  Vy^{-1}  Vz  Vz  Vy  Vz^{-1}  Vy  Vz^{-1}  Vx^{-1}  Vx^{-1} $$
 $$Vy  Vz  Vx^{-1}  Vx^{-1}  Vy  Vz^{-1}  Vx^{-1}  Vy^{-1}  Vy^{-1}  Vx^{-1}  Vy  Vz^{-1}  Vy^{-1}  Vz^{-1}  Vx  Vx  Vy  Vz  Vx^{-1}  Vy^{-1}$$
 $$Vz^{-1}  Vx  Vz^{-1}  Vy^{-1}  Vy^{-1}  Vx^{-1}  Vy^{-1}  Vy^{-1}  Vy^{-1}  Vz  Vy  Vx^{-1}  Vz  Vx^{-1}  Vy^{-1}  Vy^{-1}  Vx^{-1}  Vz$$
 $$Vx^{-1}  Vz^{-1}  Vz^{-1}  Vy  Vy  Vy  Vx^{-1}  Vy  Vy  Vy  Vz  Vy  Vx^{-1}  Vy^{-1}  Vx^{-1}  Vy^{-1}  Vz  Vz  Vz  Vy^{-1}  Vy^{-1}  Vz$$
 $$Vy  Vz  Vy^{-1}  Vx  Vx  Vx  Vy^{-1}  Vz  Vz  Vz  Vy  Vz^{-1}  Vy^{-1}  Vy^{-1}  Vy^{-1}  Vx^{-1}  Vz^{-1}  Vx^{-1}  Vz^{-1}  Vx  Vz^{-1}$$
 $$Vy^{-1}  Vx^{-1}  Vz  Vy  Vx^{-1}  Vz^{-1}  Vy^{-1}  Vx^{-1}  Vy  Vx  Vx  Vz  Vx  Vz^{-1}  Vx  Vz  Vy^{-1}  Vz  Vx^{-1}  Vy  Vz^{-1}$$
 $$Vz^{-1}  Vx  Vz^{-1}  Vx^{-1}  Vz^{-1}  Vx^{-1}  Vz  Vx  Vz^{-1}  Vx^{-1}  Vz  Vy  Vz  Vz  Vy  Vx  Vx  Vy^{-1}  Vx^{-1}  Vz^{-1}  Vx  Vy$$
 $$Vz^{-1}  Vz^{-1}  Vy  Vz^{-1}  Vy^{-1}  Vx^{-1}  Vz  Vy^{-1}  Vz^{-1}  Vy  Vx^{-1}  Vx^{-1}  Vy^{-1}  Vy^{-1}  Vy^{-1}  Vx  Vx  Vz^{-1}  $$
 $$Vx^{-1}  Vy^{-1}  Vx  Vy  Vx  Vy  Vx^{-1}  Vy  Vx^{-1}  Vx^{-1}  Vz^{-1}  Vx  Vz  Vy^{-1}  Vx^{-1}  Vy^{-1}  Vx  Vy  Vz^{-1}  Vz^{-1}  Vx $$
That is a path of size 432. The first candidate that our algorithm gives up to factor has 430 letters.  It could be a potential path but this means that the distance is optimal up to two letters. We note that the trivial lower bound for a typical element is
$$3\log_5(q)=428.5.$$ 

\subsection{Lower bound on the diameter of LPS Ramanujan graphs }
Let $a=0$, $b=1,$ which is associated to the matrix $\begin{bmatrix}i &0 \\0&-i   \end{bmatrix}$. By our correspondence in Section~\ref{quant}, the lattice point associated to this vertex is in the cups neighborhood for every $q$. In fact, this lattice point has the highest imaginary part among all the other co-volume $q$ lattice point.   Let
\begin{equation*}
\begin{split}
q=65135167346000357183003272112509282371782817584944173575600868284168\\
63929270451437126021949850746381
\end{split}
\end{equation*} The length of the shortest path from the identity to $\begin{bmatrix}i &0 \\0&-i   \end{bmatrix}$ is 571. Note that 
$4\log_5(q)=571.20,$ and recall that 
$[4\log_5(q)]$ is conjectured to be the asymptotic of the diameter of this Ramanujan graph. We refer the reader to \cite[Section 4]{Naser} for further discussion and  more numerical results.

\bibliographystyle{alpha}
\bibliography{Navigating}

\begin{thebibliography}{{Sar}15a}

\bibitem[AKS87]{Ajtai}
M.~Ajtai, J.~Komlos, and E.~Szemeredi.
\newblock Deterministic simulation in logspace.
\newblock In {\em Proceedings of the Nineteenth Annual ACM Symposium on Theory
  of Computing}, STOC '87, pages 132--140, New York, NY, USA, 1987. ACM.

\bibitem[AKS04]{primaty}
Manindra Agrawal, Neeraj Kayal, and Nitin Saxena.
\newblock Primes is in p.
\newblock {\em Annals of Mathematics}, 160(2):781--793, 2004.

\bibitem[AMM77]{Adlemann}
Leonard Adleman, Kenneth Manders, and Gary Miller.
\newblock On taking roots in finite fields.
\newblock pages 175--178, 1977.

\bibitem[BKS17]{Browning}
T.D. Browning, V.~Vinay Kumaraswamy, and R.S. Steiner.
\newblock Twisted linnik implies optimal covering exponent for $s^3$.
\newblock {\em International Mathematics Research Notices}, page rnx116, 2017.

\bibitem[EMV13]{Ellenberg}
Jordan~S. Ellenberg, Philippe Michel, and Akshay Venkatesh.
\newblock Linnik's ergodic method and the distribution of integer points on
  spheres.
\newblock In {\em Automorphic representations and {$L$}-functions}, volume~22
  of {\em Tata Inst. Fundam. Res. Stud. Math.}, pages 119--185. Tata Inst.
  Fund. Res., Mumbai, 2013.

\bibitem[HLW06]{Hoory}
Shlomo Hoory, Nathan Linial, and Avi Wigderson.
\newblock Expander graphs and their applications.
\newblock {\em Bull. Amer. Math. Soc. (N.S.)}, 43(4):439--561, 2006.

\bibitem[Len83]{Lenstra}
H.~W. Lenstra, Jr.
\newblock Integer programming with a fixed number of variables.
\newblock {\em Math. Oper. Res.}, 8(4):538--548, 1983.

\bibitem[LPS88]{Lubotzky1988}
A.~Lubotzky, R.~Phillips, and P.~Sarnak.
\newblock Ramanujan graphs.
\newblock {\em Combinatorica}, 8(3):261--277, 1988.

\bibitem[Mar88]{Margulis}
G.~A. Margulis.
\newblock Explicit group-theoretic constructions of combinatorial schemes and
  their applications in the construction of expanders and concentrators.
\newblock {\em Problemy Peredachi Informatsii}, 24(1):51--60, 1988.

\bibitem[PLQ08]{Petit2008}
Christophe Petit, Kristin Lauter, and Jean-Jacques Quisquater.
\newblock {\em Full Cryptanalysis of LPS and Morgenstern Hash Functions}, pages
  263--277.
\newblock Springer Berlin Heidelberg, Berlin, Heidelberg, 2008.

\bibitem[PS18]{Sarnak3}
Ori Parzanchevski and Peter Sarnak.
\newblock Super-golden-gates for {$PU(2)$}.
\newblock {\em Adv. Math.}, 327:869--901, 2018.

\bibitem[RS85]{Rabin}
J.~O. Rabin and Jeffrey Shallit.
\newblock Randomized algorithms in number theory.
\newblock Technical report, Chicago, IL, USA, 1985.

\bibitem[RS16]{Selinger}
Neil~J. Ross and Peter Selinger.
\newblock Optimal ancilla-free clifford+t approximation of z-rotations.
\newblock {\em Quantum Info. Comput.}, 16(11-12):901--953, September 2016.
\newblock \url{https://www.mathstat.dal.ca/~selinger/newsynth/}.

\bibitem[RS17]{Rivin}
Igor Rivin and Naser Sardari.
\newblock {\em {Quantum Chaos on random Cayley graphs}}, 2017.

\bibitem[Sar90]{Peter}
Peter Sarnak.
\newblock {\em Some applications of modular forms}, volume~99 of {\em Cambridge
  Tracts in Mathematics}.
\newblock Cambridge University Press, Cambridge, 1990.

\bibitem[{Sar}15a]{Naser2}
N.~T {Sardari}.
\newblock {Optimal strong approximation for quadratic forms}.
\newblock {\em ArXiv e-prints}, October 2015.

\bibitem[Sar15b]{Sarnak31}
Peter Sarnak.
\newblock {\em {Letter to Scott Aaronson and Andy Pollington on the
  Solovay-Kitaev Theorem}}, February 2015.
\newblock \url{https://publications.ias.edu/sarnak/paper/2637}.

\bibitem[Sar18]{Naser}
Naser~T. Sardari.
\newblock Diameter of ramanujan graphs and random cayley graphs.
\newblock {\em Combinatorica}, Aug 2018.

\bibitem[Sch85]{Schoof}
Ren\'e Schoof.
\newblock Elliptic curves over finite fields and the computation of square
  roots mod {$p$}.
\newblock {\em Math. Comp.}, 44(170):483--494, 1985.

\bibitem[Sha73]{Shanks}
Daniel Shanks.
\newblock Five number-theoretic algorithms.
\newblock pages 51--70. Congressus Numerantium, No. VII, 1973.

\bibitem[Sho92]{Shoup}
Victor Shoup.
\newblock Searching for primitive roots in finite fields.
\newblock {\em Math. Comp.}, 58(197):369--380, 1992.

\end{thebibliography}

\end{document}